\documentclass[letterpaper, 10 pt, conference]{ieeeconf}  

\IEEEoverridecommandlockouts                              
\usepackage{graphics} 
\usepackage{epsfig} 
\usepackage{times} 
\usepackage{amsmath} 
\usepackage{amssymb}  

\usepackage{amsfonts}
\usepackage[caption=false,font=footnotesize]{subfig} 
\usepackage{mathtools}
\usepackage{algorithm}
\usepackage{algpseudocode}
\usepackage{xcolor}

\newtheorem{assumption}{Assumption}

\newtheorem{remark}{Remark}

\title{\LARGE \bf
Differential Dynamic Programming for the Optimal Control Problem\\ with an Ellipsoidal Target Set and Its Statistical Inference
}

\author{Sungjun Eom and Gyunghoon Park$^*$
\thanks{Sungjun Eom and Gyunghoon Park are with the School of Electrical and Computer Engineering,
        University of Seoul, Republic of Korea
        {\tt\small \{junjun607,gyunghoon.park\}@uos.ac.kr}}%
\thanks{This work was supported by the Institute of Information \& Communications
  Technology Planning \& Evaluation (IITP) grant funded by the Korea government (MSIT)
  (No. RS-2025-02305523).}%
}

\begin{document}

\maketitle
\thispagestyle{empty}
\pagestyle{empty}


\begin{abstract}
This work addresses an extended class of optimal control problems where a target for a system state has the form of an ellipsoid rather than a fixed, single point. 
As a computationally affordable method for resolving the extended problem, we present a revised version of the differential dynamic programming (DDP), termed the differential dynamic programming with ellipsoidal target set (ETS-DDP).  
To this end, the problem with an ellipsoidal target set is reformulated into an equivalent form with the orthogonal projection operator, yielding that the resulting cost functions turn out to be discontinuous at some points.
As the DDP usually requires the differentiability of cost functions, in the ETS-DDP formulation we locally approximate the (nonsmooth) cost functions to smoothed ones near the path generated at the previous iteration, by utilizing the explicit form of the orthogonal projection operator. 
Moreover, a statistical inference method is also presented for designing the ellipsoidal target set, based on data on admissible target points collected by expert demonstrations.
Via a simulation on autonomous parking of a vehicle, it is seen that the proposed ETS-DDP efficiently derives an admissible state trajectory while running much faster than the point-targeted DDP, at the expense of optimality.
\end{abstract}


\section{Introduction}

Among several promising methodologies for trajectory optimization and/or optimal control, the differential dynamic programming (DDP) has been particularly renowned in these days, as a computationally affordable variant of the well-known dynamic programming (DP). 
Unlike the conventional DP that searches all possibilities in a single step, the DDP recursively computes the optimal trajectory with the help of the local approximation of the Bellman's equation. 

Since the birth of the DDP in 1980s, a tremendous number of research efforts have been made in order to extend the class of problems to be resolved via the DDP and its variants. 
The DDP was initially developed for unconstrained  problems \cite{mayne1966second}, and then theoretical guarantees on its convergence were explored in \cite{shoemaker1990proof, liao2002convergence}. 
After a while, a number of variants have been proposed to deal with constrained problems, with box-type input constraints \cite{tassa2014control} and more recently, with linear or nonlinear constraints  \cite{pavlov2021interior, mastalli2022feasibility, almubarak2022safety}.
On the other hand, some researchers have attempted to incorporate the DDP algorithms into more complicated scenarios, such as distributed and/or decentralized settings \cite{saravanos2023distributed, huang2023decentralized}.

Despite such huge success in the literature, the existing DDP algorithms are still limited in their usage, to the cases when the target of the control task is given as a single state. 
In real-world scenarios such as parking a vehicle, however, the task is often defined in a more abstract manner, e.g., driving the terminal state to enter a certain set (not reach a point). 
From this perspective, it would be more realistic to consider the target in the optimal control problem as a set of acceptable states, rather than a single point. 

\begin{figure}[htbp!]
  \centering
  \subfloat[With a point-valued target]{%
    \includegraphics[width=0.5\linewidth]{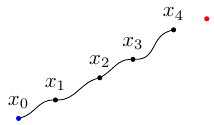}}
  \hfill
  \subfloat[With a set-valued target]{%
    \includegraphics[width=0.5\linewidth]{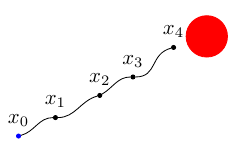}}
  \caption{Conceptual illustration of control problem with point-valued and set-valued targets: $x_t$ represents the system state, and the target is marked in red.}
  \label{fig:ocp}
\end{figure}

As an effective solver for the optimal control problem with {\it {a} target set}, in this paper we propose a new formulation of the DDP, termed \textbf{D}ifferential \textbf{D}ynamic \textbf{P}rogramming with \textbf{E}llipsoidal \textbf{T}arget \textbf{S}et (ETS-DDP). 
The core idea of the ETS-DDP is twofold. 
Provided that the target set has the form of an ellipsoid (that is inherently convex), it is first seen that the optimal control problem with the target set can be equivalently represented into more DDP-friendly form, with the help of the orthogonal projection operator. 
While orthogonal projection introduces the nonsmoothness of the cost in the equivalent problem, we tackle this issue by approximating the nonsmooth part in a local region by utilizing the path from the previous iteration in the backward pass, by which the ETS-DDP can be derived. 
In addition to the structural modification of the DDP algorithm, we further present a statistical way to construct the ellipsoidal target set based on pre-collected data from expert demonstrations.  
It is observed via simulation on autonomous parking that, compared with the point target-based conventional DDPs, the proposed framework on the ETS-DDP offers an opportunity to accelerate the DDP algorithm as well as reduce human intervention during decision process, at the expense of losing optimality. 

\section{Preliminaries on Differential Dynamic Programming}

We begin by summarizing how the conventional differential dynamic programming (DDP) works. 
The DDP basically aims to solve the following optimal control problem defined for given initial state $x_0\in {\mathbb{R}}^{n}$:
\begin{subequations}\label{eq:ocp}
    \begin{align}
    \min_{\mathbf{u}}~&~ \left\{\sum_{t=0}^{T-1} L(x_t-c, u_t) + \phi(x_T-c)\right\}, \\
    {\rm s.t.}~&~ x_{t+1} = f(x_t, u_t),\quad \forall t=0,\dots,T-1
\end{align}
\end{subequations}
where $x_t\in {\mathbb{R}}^n$ and $u_t\in {\mathbb{R}}^l$ are the state and input of the dynamical system $x_{t+1} = f(x_t,u_t)$, $T$ is the horizon length, $c\in {\mathbb{R}}^n$ is the target for the state, ${\bf u}:=[u_0;\cdots;u_{T-1}]\in {\mathbb{R}}^{Tl}$ is the stacked vector of $u_t$. 
In addition, positive semi-definite functions $L:{\mathbb{R}}^n\times {\mathbb{R}}^l\rightarrow{\mathbb{R}}$ and $\phi:{\mathbb{R}}^n\rightarrow {\mathbb{R}}$ represent the stage cost and terminal cost functions, respectively. 
It has been assumed in the relevant literature that, the functions $L$, $\phi$ and $f$ are sufficiently smooth, and the target $c$ for the state is chosen as a fixed point. 
For ease of explanation, let
\begin{align*}
    L^c (x,u):= L(x-c,u),\quad \phi^c (x):= \phi(x-c). 
\end{align*}

For the (point-targeted) optimal control problem \eqref{eq:ocp} and the corresponding value function $V(x_t)$, the DDP algorithm aims to find its optimal solution by iteratively solving to the {\it Bellman's equation}
\begin{align}
    V(x_t) &= \min_{u_t}\{L^c(x_t, u_t) + V(x_{t+1})\}\notag\\
    &=\min_{u_t}\{L^c(x_t, u_t) + V(f(x_t, u_t)\}\label{eq:V}
\end{align}
backward in time. 
In \eqref{eq:V}, the function to be minimized
\begin{equation*}
    Q(x_t,u_t):=L^c(x_t,u_t)+V(f(x_t,u_t))
\end{equation*}
is usually called the {\it Q-function}.
At each iteration, the DDP updates from $u_t$ to $u_t+\delta u_t$ in a way that the resulting Q-function is minimized over a local region, similar to Gauss-Newton methods. 
To proceed, we approximate $Q(x,u)$ near $(x,u)=(x_t,u_t)$ by the second-order Taylor expansion as 
\begin{equation*} 
\begin{split}
    & Q(x+\delta x, u+\delta u)\\
    & \approx Q(x,u) + Q_x^\top \delta x+ Q_u^\top\delta u\\
    & \quad +\delta x^\top Q_{xu}\delta u + \frac{1}{2}\delta x^\top Q_{xx}\delta x + \frac{1}{2}\delta u^\top Q_{uu} \delta u\\
    &=:\widehat{Q}(\delta x, \delta u),
\end{split}
\end{equation*}
where $Q_{x}$, $Q_u$, $Q_{ux}$, $Q_{xx}$, $Q_{uu}$ are partial derivatives of $Q$, satisfying
\begin{subequations}\label{eq:DDP-Q}
	\begin{align}
		Q_x&=L_x^c+f_x^\top V_x(x_{t+1}),\\
		Q_u&=L_u^c+f_u^\top V_x(x_{t+1}),\\
		Q_{ux}&=L_{ux}^c+f_{ux}\cdot V_x(x_{t+1})+f_u^\top V_{xx}(x_{t+1})f_x,\\
		Q_{xx}&=L_{xx}^c+f_{xx}\cdot V_x(x_{t+1})+f_x^\top V_{xx}(x_{t+1})f_x,\\
		Q_{uu}&=L_{uu}^c+f_{uu}\cdot V_x(x_{t+1})+f_u^\top V_{xx}(x_{t+1})f_u. 
	\end{align}
\end{subequations}
It is noted that, in the above, the second-order derivatives $f_{xx}$, $f_{ux}$, $f_{uu}$ of the vector-valued function $f$ are in fact third-order tensors, and the dot $\cdot$ denotes a {\it tensor dot product}.

After some computations, one can readily see that
\begin{align}
	\delta u = \delta u^* := \underbrace{-Q_{uu}^{-1}Q_u}_{=:k}+ \underbrace{(-Q_{uu}^{-1}Q_{ux})}_{=:K}\delta x \label{eq:u*}
\end{align}
minimizes the approximated Q-function $\widehat{Q}(\delta x,\delta u)$

Similarly, we can also perform a Taylor expansion of the value function $V$ up to the second order:
\begin{align*}
    V(x+\delta x)&\approx V(x)+ V_x^\top \delta x+\frac{1}{2}\delta x^\top V_{xx}\delta x
\end{align*}
where the partial derivatives $V_x$ and $V_{xx}$ of the value function $V$ are related to those of $Q$ as
\begin{subequations}\label{eq:DDP-VQ}
	\begin{align}
 	 V_x&=Q_x-Q_{xu}Q_{uu}^{-1}Q_u,\\
		V_{xx}&=Q_{xx}-Q_{xu}Q_{uu}^{-1}Q_{ux}
	\end{align}
\end{subequations}
which are utilized at the next iteration for computing \eqref{eq:DDP-Q}.

Summarizing so far, a psuedocode for the conventional DDP is presented in Algorithm~\ref{al:DDP} for the unconstrained problem \eqref{eq:ocp}. 
Note that one iteration of the DDP consists of a backward pass and a forward pass. 
In the backward pass, the DDP recursively computes the partial derivatives \eqref{eq:DDP-Q} of $Q$ backward in time, with which the feedback gain $K$ and the feedforward term $k$ in \eqref{eq:u*} is computed.
Once the process for the backward pass ends up, the forward pass begins to update the system state $x_t$ forward in time (i.e., from $t=0$ to $t=T$) by utilizing the system model $x_{t+1}=f(x_t,u_t)$ as well as the locally-optimal input $u_t + \delta u_t$.
More details can be found in the relevant works.

\begin{algorithm}[t]
	\caption{Pseudocode for conventional DDP (for the unconstrained problem \eqref{eq:ocp})}\label{al:DDP}
	\renewcommand{\algorithmicrequire}{\textbf{Input:}}
	\begin{algorithmic}[1]
		\Require{horizon length \(T\), vector field $f(\cdot)$ of system, initial state $x_0$, stage cost $L(\cdot)$, terminal cost $\phi(\cdot)$}
            \State Initialize $\{x_t\}$ and$\{u_t\}$ by zeros.
		\For {$m = 0,1,2,\ldots$}\textcolor{gray}{\Comment{iteration}}
		\State compute $V_x, V_{xx}$ at $x_T$
		\For {$t = T-1,T-2,\ldots,0$} \textcolor{gray}{\Comment{backward pass}}
		\State compute Q derivatives using \eqref{eq:DDP-Q} at $x_t$
		\State compute $k_t$ and $K_t$ using \eqref{eq:DDP-Q}
		\State compute $V_{x,t}$ and $V_{xx,t}$ for the next step by \eqref{eq:DDP-VQ}
		\EndFor
		\State $u_0\leftarrow u_0+k_0$
		\For {$t = 0,1,2,\ldots,T-1$}\textcolor{gray}{\Comment{forward pass}}
		\State update state $x_{t+1}$ using dynamical relationship $x_{t+1}=f(x_t,u_t)$
		\State update $u_{t+1}$ using \eqref{eq:u*}
		\EndFor
		\If {the amount of change of the total cost is smaller than the threshold}
		escape the loop
		\EndIf
		\EndFor
	\end{algorithmic}
\end{algorithm}

\section{Differential Dynamic Programming with Ellipsoidal Target Set}

\subsection{Optimal Control Problem with Target Set}

In this paper, we address an extended form of the optimal control problem \eqref{eq:ocp}, in the sense that the target $c$ for the control task is no longer a fixed point anymore: 
instead, the target is supposed to be a subset of the state space ${\mathbb{R}}^n$, denoted by $\mathcal{C}$. 
With additional {\it decision variables} $c_t$, $t=0,\dots,T$, the optimal control problem of the form \eqref{eq:ocp} can be relaxed into 
\begin{subequations}\label{eq:ocp2}
    \begin{align}
    \min_{\mathbf{u}, {\bf c}}~&~ \left\{\sum_{t=0}^{T-1} L(x_t-c_t, u_t) + \phi(x_T-c_T)\right\}, \!\!\\
    {\rm s.t.}~&~ x_{t+1} = f(x_t, u_t),\\
    &~ c_t \in \mathcal{C}
\end{align}
\end{subequations}
where ${\bf c}:=[c_0,\dots,c_T]$ consists of the extra variables $c_t$.
Throughout this paper, we call the problem in \eqref{eq:ocp2} an {\it optimal control problem with target set $\mathcal{C}$}. 
In the rest of this section, we will propose a variant of the DDP, with which the solution of the relaxed problem \eqref{eq:ocp2} can be numerically derived. 
It is worthwhile noting that applying the conventional DDP (i.e., Algorithm~\ref{al:DDP}) into \eqref{eq:ocp2} is not straightforward, mainly because the new variable $c_t$ does not have a generating dynamics, which cannot be handled in the DDP formalism directly. 

Before closing this section, it is pointed out that such a relaxation in the cost function has a merit for a particular class of industrial applications, where a control task to be satisfied is defined by a set (e.g., $x_T\in \mathcal{C}$).
One of the remarkable examples is the autonomous parking of a vehicle. 
Indeed, one can determine the success of the parking task by checking whether or not the state of a vehicle is located inside a set (that represents a parking space), in which it is not that important how close the state of the vehicle is with a specific fixed point. 
We will revisit this example later with computer-aided simulations.

\subsection{Differential Dynamic Programming with Ellipsoidal Target Set: ETS-DDP}

As this work is the first attempt in the literature to develop a DDP for the OCP \eqref{eq:ocp2} with a target set, we present an assumption on $\mathcal{C}$, which makes the considered problem relatively simple. 

\begin{assumption}
    The target set $\mathcal{C}$ in \eqref{eq:ocp2} is an ellipsoid of the form 
\begin{equation}\label{eq:D}
    {\mathcal{C}}= \left\{c \in {{\mathbb{R}}^n} :\sqrt{(c-o)^\top \Sigma^{-1} (c-o)}< r \right\}
\end{equation}
    where $\Sigma= \Sigma^\top >0$, $r>0$ is the radius, and $o$ is the center of the ellipsoid.  $\hfill\square$ 
\end{assumption}

It is noted in advance that a synthesis method for such $\mathcal{C}$ will be discussed shortly. 
On the other hand, the following assumption requires a mild condition for the cost functions in \eqref{eq:ocp2}. 
\begin{assumption}\label{asm:cost}
    For any $u\in {\mathbb{R}}^l$ and $a,b\in {\mathbb{R}}^n$ satisfying $0\leq  \|a\|\leq \|b\|$, 
    \begin{align*}
        L(a,u)\leq L(b,u),\quad \phi(a)\leq \phi(b).
    \end{align*}
    $\hfill\square$
\end{assumption}

We are now ready to present the core idea of reducing the number of decision variables in \eqref{eq:ocp2}. 
By Assumption~\ref{asm:cost}, for any $c\in \mathcal{C}$ and $u\in {\mathbb{R}}^n$, the cost functions in \eqref{eq:ocp2} have a lower bound
\begin{subequations}
\begin{align}
    L(x-c,u) & \geq L\big(x-P_{\mathcal{C}}(x),u\big),\\
    \phi(x-c) & \geq \phi \big(x - P_{\mathcal{C}}(x)\big)
\end{align}
\end{subequations}
where $P_{\mathcal{C}}:{\mathbb{R}}^n \rightarrow {\mathcal{C}}$ is the {\it orthogonal projection operator} with respect to $\mathcal{C}$, defined as \cite{beck2014introduction}:
\begin{equation*} 
    P_{\mathcal{C}}(x)=\arg\min\{||y-x||^2: y\in {\mathcal{C}}\}.
\end{equation*}
Note that the value $P_{\mathcal{C}}(x)$ is the closest point in ${\mathcal{C}}$ to $x$, and is uniquely determined if $\mathcal{C}$ is an ellipsoid as in \eqref{eq:D} (and thus convex). 
Then it is readily concluded that, the optimal control problem \eqref{eq:ocp2} is in fact equivalent to
\begin{subequations}\label{eq:newocp}
    \begin{align}
    \min_{\mathbf{u}}~&~ \left\{\sum_{t=0}^{T-1} L(x_t-P_{\mathcal{C}}(x_t), u_t) + \phi(x_T-P_{\mathcal{C}}(x_T))\right\}, \!\!\\
    {\rm s.t.}~&~ x_{t+1} = f(x_t, u_t),\quad \forall t=0,\dots,T-1
\end{align}
\end{subequations}
where the decision variables reduce from $(u_t,c_t)$ to $u_t$, which is in fact the same as in the original (point-targeted) problem \eqref{eq:ocp}. 
In what follows, let
\begin{subequations}\label{eq:LC_phiC}
\begin{align}
	L^{\mathcal{C}}(x,u)  & := L(x-P_{\mathcal{C}}(x), u),\\
	\phi^{\mathcal{C}}(x) & := \phi(x -P_{\mathcal{C}}(x)).
\end{align}
\end{subequations}

When it comes to implementation of the DDP for the equivalent problem \eqref{eq:newocp} as in Algorithm~\ref{al:DDP}, one may face another problem that $L^{\mathcal{C}}$ and $\phi^{\mathcal{C}}$ are not differentiable at $x\in \partial \mathcal{C}$. 
With this kept in mind, in the following we propose a variant of the DDP termed {\it DDP with ellipsoidal target set} (ETS-DDP), by locally approximating $L^{\mathcal{C}}$ and $\phi^{\mathcal{C}}$ into differentiable forms based on the path computed in the previous iteration of the DDP. 
To this end, we express the result of the orthogonal projection operator as
\begin{equation*}
    P_{\mathcal{C}}(x)=\begin{cases}
    x-o, & \forall x \in \mathcal{C},\\ 
    \Sigma^{1/2}(\Sigma+\lambda_K I)^{-1}\Sigma^{1/2}(x-o)+o, 
    & \forall x \notin \mathcal{C}
    \end{cases}
\end{equation*}
where the ellipsoidal structure of $\mathcal{C}$ has been used. $\lambda_K$ is a Lagrange multiplier at its optimum. Full derivation of the orthogonal projection operator when $x\notin \mathcal C$ is given in Appendix. If we define $y(x):=\Sigma^{1/2}(\Sigma+\lambda_K I)^{-1}\Sigma^{1/2}(x-o)+o$,
then near a given $x^-$ (that implicitly represents $x_t$ computed in the previous iteration), possibly nonsmooth cost functions $L^{\mathcal{C}}(x,u)$ and  $\phi^{\mathcal{C}}(x)$ in \eqref{eq:LC_phiC} are locally approximated as
\begin{subequations}\label{eq:hatLC}
    \begin{align}
    & \widehat{L}^{\mathcal{C}}(x,u;x^-)\\
    & := \begin{cases}
    L(o,u), & \forall x^- \in \mathcal{C},\\
    L\left(x-y(x),u\right), &\forall x^- \notin \mathcal{C},
    \end{cases}\notag\\
    & \widehat{\phi}^{\mathcal{C}}(x;x^-)\\
    & := \begin{cases}
    \phi(o), &\forall x^- \in \mathcal{C},\\
    \phi\left(x-y(x)\right), &\forall x^- \notin \mathcal{C},
    \end{cases}\notag
    \end{align}
\end{subequations}
respectively. 
In view of the DDP algorithm, as $L^c$ is replaced with the smoothed function $\widehat{L}^{\mathcal{C}}$, the  formula \eqref{eq:DDP-Q} for computing the partial derivatives of $Q(x,u)$ must be modified as:
\begin{subequations}\label{eq:ETS-DDP_Q}
\begin{align}
    Q_x & = \widehat{L}^{\mathcal{C}}_x+f_x^\top V_x(x_{t+1}),\!\\
    Q_u & = \widehat{L}^{\mathcal{C}}_u+f_u^\top V_x(x_{t+1}),\!\\
    Q_{ux}&=\widehat{L}^{\mathcal{C}}_{ux}+f_{ux}\cdot V_x(x_{t+1})+f_u^\top V_{xx}(x_{t+1})f_x,\!\!\!\\
    Q_{xx}&=\widehat{L}^{\mathcal{C}}_{xx}+f_{xx}\cdot V_x(x_{t+1})+f_x^\top V_{xx}(x_{t+1})f_x,\!\!\!\\
    Q_{uu}&=\widehat{L}^{\mathcal{C}}_{uu}+f_{uu}\cdot V_x(x_{t+1})+f_u^\top V_{xx}(x_{t+1})f_u,  \!\!\!
\end{align}
\end{subequations}
where $\widehat{L}^{\mathcal{C}}_x$, $\widehat{L}^{\mathcal{C}}_u$, $\widehat{L}^{\mathcal{C}}_{ux}$, $\widehat{L}^{\mathcal{C}}_{xx}$, and  $\widehat{L}^{\mathcal{C}}_{uu}$ are partial derivatives of the smoothed functions \eqref{eq:hatLC}.

In summary, for the optimal control problem \eqref{eq:ocp2} with ellipsoidal target set $\mathcal{C}$ in \eqref{eq:D}, the ETS-DDP basically runs similarly to the conventional DDP (e.g., Algorithm~\ref{al:DDP}), but $L^c$, $\phi$, and partial derivatives \eqref{eq:DDP-Q} being properly replaced by $\widehat{L}^{\mathcal{C}}$ and $\widehat{\phi}^{\mathcal{C}}$ in \eqref{eq:hatLC}, and \eqref{eq:ETS-DDP_Q}. 

{
	\begin{remark}
		We have so far discussed that for extending the class of the target (from a point $c$ to a set $\mathcal{C}$), it is enough to replace target-related functions with their modified versions, while the structure of the algorithm is mostly maintained.   
		Therefore, it would be natural to apply the core idea of the ETS-DDP to not only unconstrained problem \eqref{eq:newocp} but also constrained ones, for example, having extra box-type input constraints $\underline{u}\leq u_t \leq \overline{u}$.
		It will be seen in the simulation part that application to constrained case is indeed possible. $\hfill\square$  
	\end{remark}
}

\subsection{Synthesis of an Ellipsoidal Target Set via Expert Demonstrations}\label{subsec:ellipsoid}

As the ETS-DDP proposed in the previous subsection requires the ellipsoidal target set $\mathcal{C}$ to be pre-determined, it is also necessary for completeness of the overall framework to provide a proper design method for $\mathcal{C}$.
In this subsection, we present a {\it statistical approach} to synthesizing $\mathcal{C}$ based on a set of experiments demonstrated by an expert, expected to reduce human intervention in the end.  

To this end, it is first assumed that we gather $N$ independent and identically distributed random data points $c^{(i)}$, $i=1,\dots,N$.
More precisely, each data point $c^{(i)}$ is expected to be collected by an expert as a {\it possible (terminal) state with success of control}  in a way that, if the state $x_t$ of the system eventually turns out to be $c^{(i)}$ for some $t$, then one can judge that the control task is achieved.  
Thus, if $\mathcal{C}$ is selected to contain a large number of such $c^{(i)}$ while considering inherent randomness of an expert, then it  is empirically reasonable to formulate the optimal control problem as in \eqref{eq:newocp} with $\mathcal{C}$ being a target set.

For ease of discussion to come, the following assumption is made. 
\begin{assumption}
	The data points  $c^{(i)}$ follows the multivariate normal distribution ${\mathcal{N}}(o, \Sigma)$. $\hfill\square$
\end{assumption}

A natural consequence of this presumption is that the square of the {\it Mahalanobis distance}
\begin{align*}
d_M\big(c^{(i)};o,\Sigma):=\sqrt{(c^{(i)}-o)^\top \Sigma^{-1} (c^{(i)}-o)}
\end{align*}
follows the chi-squared distribution with $n$ degrees of freedom \cite{rencher1997methods}: that is,
\begin{equation}
    \big(d_M\big(c^{(i)};o,\Sigma)\big)^2 \sim \chi^2(n).\label{eq:dM_convergence}
\end{equation}
With a number of data points $c^{(i)}$ above, we compute the sample mean $\overline{c}$ and the sample variance $S$ as 
\begin{align*}
    \bar{c} :=\frac{1}{N}\sum_{i=1}^{N}c ^{(i)},~
    S :=\frac{1}{N-1}\sum_{i=1}^{N}(c^{(i)}-\overline{c})(c^{(i)}-\overline{c})^\top.
\end{align*}
The strong law of large numbers states that, as the number $N$ of the samples is sufficiently large, the sample mean $\overline{c}$ and the sample variance matrix $S$ approximates the true values in the sense that
\begin{align*}
    \overline{c} &\xrightarrow{\mathrm{a.s.}} o \text{\quad and\quad}
    S \xrightarrow{\mathrm{a.s.}} \Sigma
\end{align*}
where $\mathrm{a.s.}$ stands for ``almost sure" convergence.
Similarly, the {\it sample} Mahalanobis distance $d_M(c^{(i)};\overline{c},\overline{\Sigma})$ also has a convergent property as
\begin{align}
	d_M(c^{(i)};\overline{c},S) \xrightarrow{\mathrm{a.s.}} d_M(c^{(i)};o,\Sigma)\label{eq:dM_convergence2}
\end{align}
where the right-hand side term satisfies \eqref{eq:dM_convergence}.  
The discussions made so far introduce a statistical design guideline for the parameters $o$, $r$, and $\Sigma$ of the ellipsoidal target set $\mathcal{C}$ in \eqref{eq:D}, based only on the data points $c^{(i)}$ collected from expert demonstrations. 
From now on, we assume the approximation \eqref{eq:dM_convergence2} holds and take $o$ and $\Sigma$ in \eqref{eq:D} as the sample mean $\bar{c}$ and sample variance $S$.

For selection of the radius $r$ that determines the volume of $\mathcal{C}$, we here bring the concept of {\it (prediction) interval estimation} on a single data $c^{(i)}$ into the picture. 
To this end, choose $\alpha$ within the range $0<\alpha<1$, which will serve as a {\it significance level}. 
With $\alpha$ given, the radius $r$ is determined as
\begin{align}\label{eq:parameter_C_2}
	r = \sqrt{\chi^2_{\alpha}(n)}
\end{align}
where for any $\beta\in (0,1)$, $\chi_\alpha^2(n)$ stands for the $(1-\alpha)$-quantile of the chi-squared distribution with $n$ degrees of freedom. 
It should be pointed out that, by definition of the $(1-\alpha)$-quantile, ${\mathbb{P}}\big( d_M\big(c^{(i)};\bar{c},S)\leq r \big) = { 1-\alpha}.$
This {\it roughly} says that the resultant ellipsoid $\mathcal{C}$ contains about $1-\alpha$ ratio of data points $c^{(i)}$.

\begin{remark}
	Basically, the specific value of $\alpha$ (and thus $r$) depends on the choice of an user. 
	The proposed guideline accounts for innate randomness from expert demonstration, implicitly yielding that one may be confident on the accuracy of the expert demonstrations with the level $\alpha$ of significance. 
	It should be emphasized here that, this in fact does not mean that the probability of  a single data $c^{(i)}$ being in $\mathcal{C}$ equals $1-\alpha$, which is a common misconception. 
	$\hfill\square$
\end{remark}

\begin{figure*}[!t]
  \centering
  \subfloat[Case of acceptance\label{fig:softwarea}]{
    \includegraphics[width=0.3\linewidth]{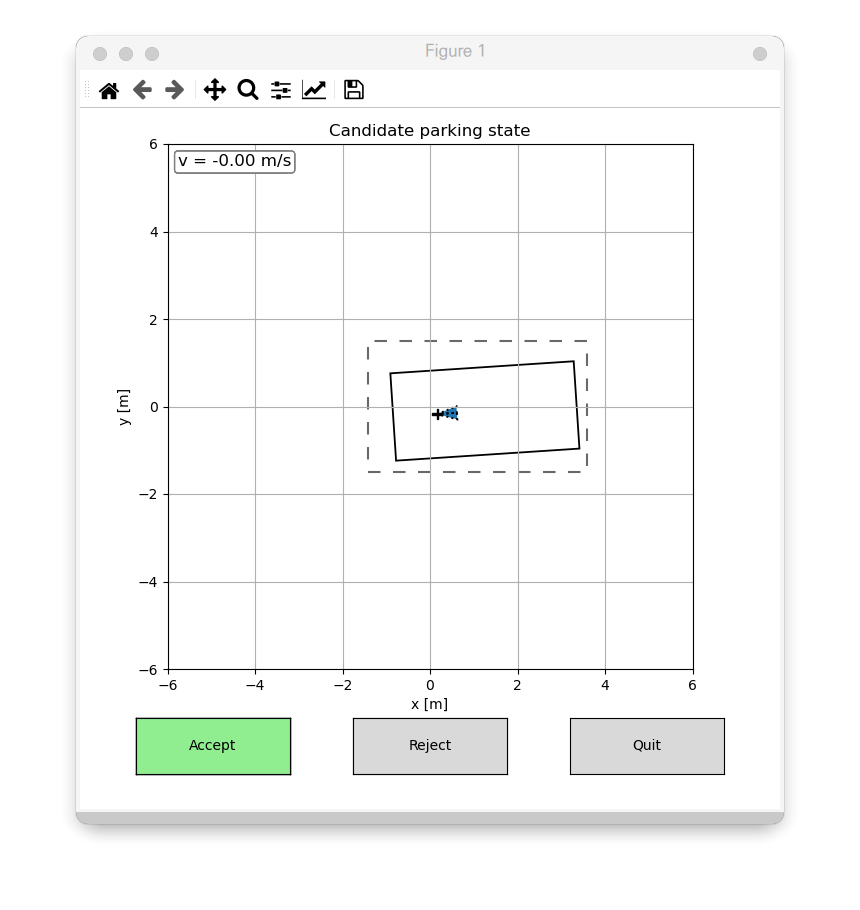}
  }\hfill
  \subfloat[Case of rejection\label{fig:softwareb}]{
    \includegraphics[width=0.3\linewidth]{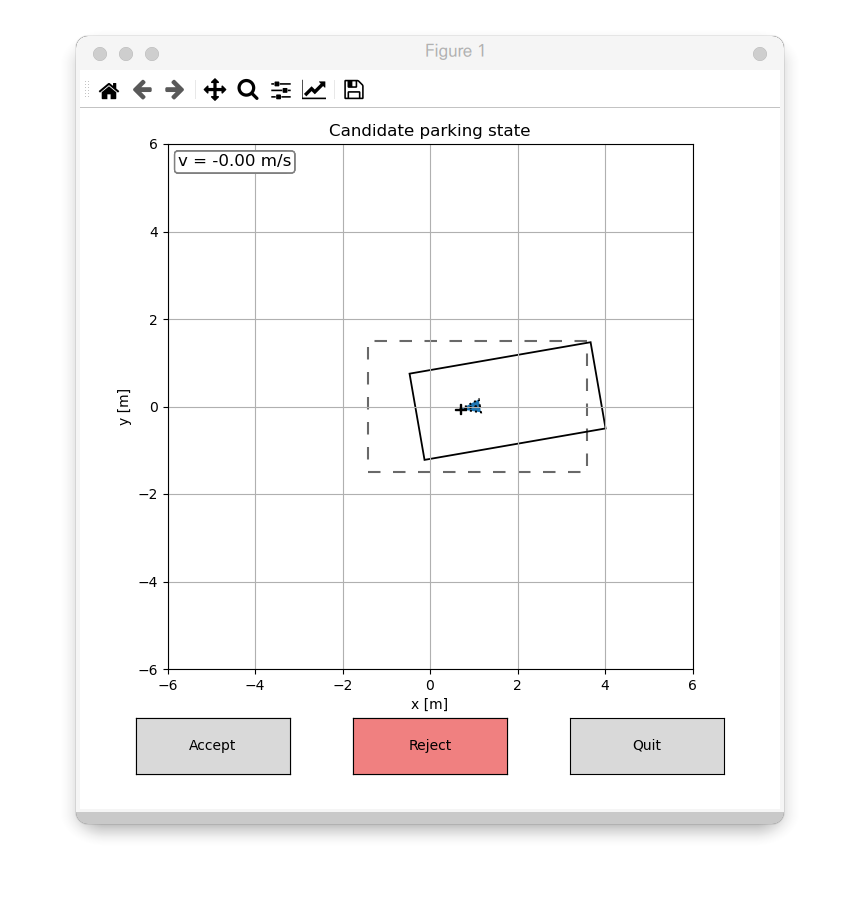}
  }\hfill
  \subfloat[Collected data $c^{(i)}$ from expert demonstration\label{fig:softwarec}]{
    \includegraphics[width=0.3\linewidth]{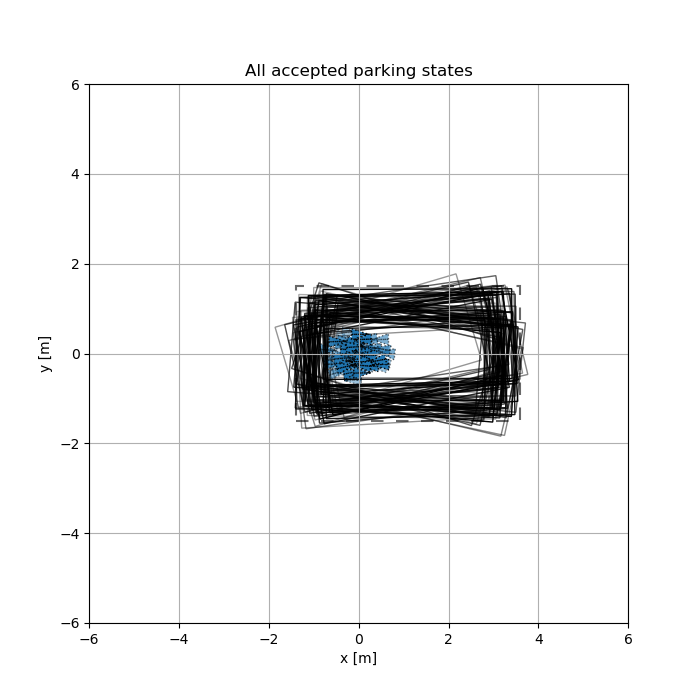}
  }
  \caption{Data collecting software.}
  \label{fig:software}
\end{figure*}

\section{Simulation} \label{sec:sim}

\subsection{Simulation Setup}
In the following, we will demonstrate our proposed method on the 2D car parking scenario \cite{tassa2014control}.
In discrete time, the kinematics model of a vehicle to be controlled is given by 
\begin{align*}
	\underbrace{\begin{bmatrix}
			p^x_{t+1}\\
			p^y_{t+1}\\
			\theta_{t+1}\\
			v_{t+1}
	\end{bmatrix}}_{=x_{t+1}}
	=
	\underbrace{
		\begin{bmatrix}
			p^x_t + b \cos (\theta_t)\\
			p^y_t + b \sin (\theta_t)\\
			\theta_t + \sin^{-1} \left(  (\Delta/d)\sin(\omega_t) v_t \right)\\
			v_t + \Delta a_t
		\end{bmatrix}
	}_{=:f(x_t,u_t)}
\end{align*}
where $x_t:=(p^x_t,p^y_t, \theta_t, v_t)\in {\mathbb{R}}^4$ stands for the state, and $u_t:=(\omega_t,a_t)\in {\mathbb{R}}^2$ is the control input.  
Here, $p^x_t$ and $p^y_t$ represent the position of the vehicle in $x$- and $y$-directions, respectively, $\theta_t$ denotes the orientation angle, $v_t$ and $\omega_t$ are the linear and  angular velocities of the vehicle, respectively, and $a_t$ is the (linear) acceleration. 
The parameters $b$ and $d$ are some constants and the detailed values can be found in  \cite{tassa2014control}. 

For a comparison, we here construct two types of DDPs.
One is the conventional DDP having the target $c$ as a fixed, single point, while the other is the proposed ETS-DDP with an ellipsoidal target set $\mathcal{C}$ (which will be explained in details in the upcoming subsection). 
Following the active set-based method presented in \cite{tassa2014control}, both DDPs are constructed in order to deal with the {\it box constraints} $\omega \in [-0.5,0.5]$ and $a \in [-2,2]$ on the control input. 
 
For simplicity, we set the target point $c$ of the conventional DDP as the origin $c=0$ (as in \cite{tassa2014control}).
In other words, the control task is to drive the state $x_t$ of the vehicle near the origin $x=c=0$ at $t=T$. 
Similarly, the stage cost $L$ and terminal cost $\phi$ in the optimal control problems \eqref{eq:ocp} and \eqref{eq:ocp2} are chosen as the same as in \cite{tassa2014control}: that is, with the {\it Huber-type function}  $h_\mu(\cdot)=\sqrt{(\cdot)^2+\mu^2}-\mu$ that smoothly approximates the absolute-valued function $|\cdot|$, we set
\begin{align*}
	L(x,u) & = {\sf q}_1 h_{\mu_1}(p^x) + {\sf q}_2 h_{\mu_2}(p^y) + {\sf r}_1 \omega^2 + {\sf r}_2 a^2,\\
	\phi(x) & = h_{\mu_1}(p^x) + h_{\mu_2}(p^y) + h_{\mu_3}(\theta) + h_{\mu_4}(v)
\end{align*}
where the weights are taken as ${\sf q}_1 = {\sf q}_2 = 0.01$, ${\sf r}_1 = 0.01$, ${\sf r}_2 = 0.0001$, while the smoothness scales are given as $\mu_1=\mu_2=0.1$, $\mu_3=0.01$, and $\mu_4=1$.
In both DDPs, we set the horizon length $T$ as $500$. 

The simulation was performed with MATLAB (R2024b) on Red Hat Enterprise Linux 9.5 and 13th Gen Intel Core i9-13900K (5.80 GHz) CPU. 

\begin{figure}[!t]
  \centering
  \subfloat[Top view\label{fig:sima}]{
    \includegraphics[width=0.8\linewidth]{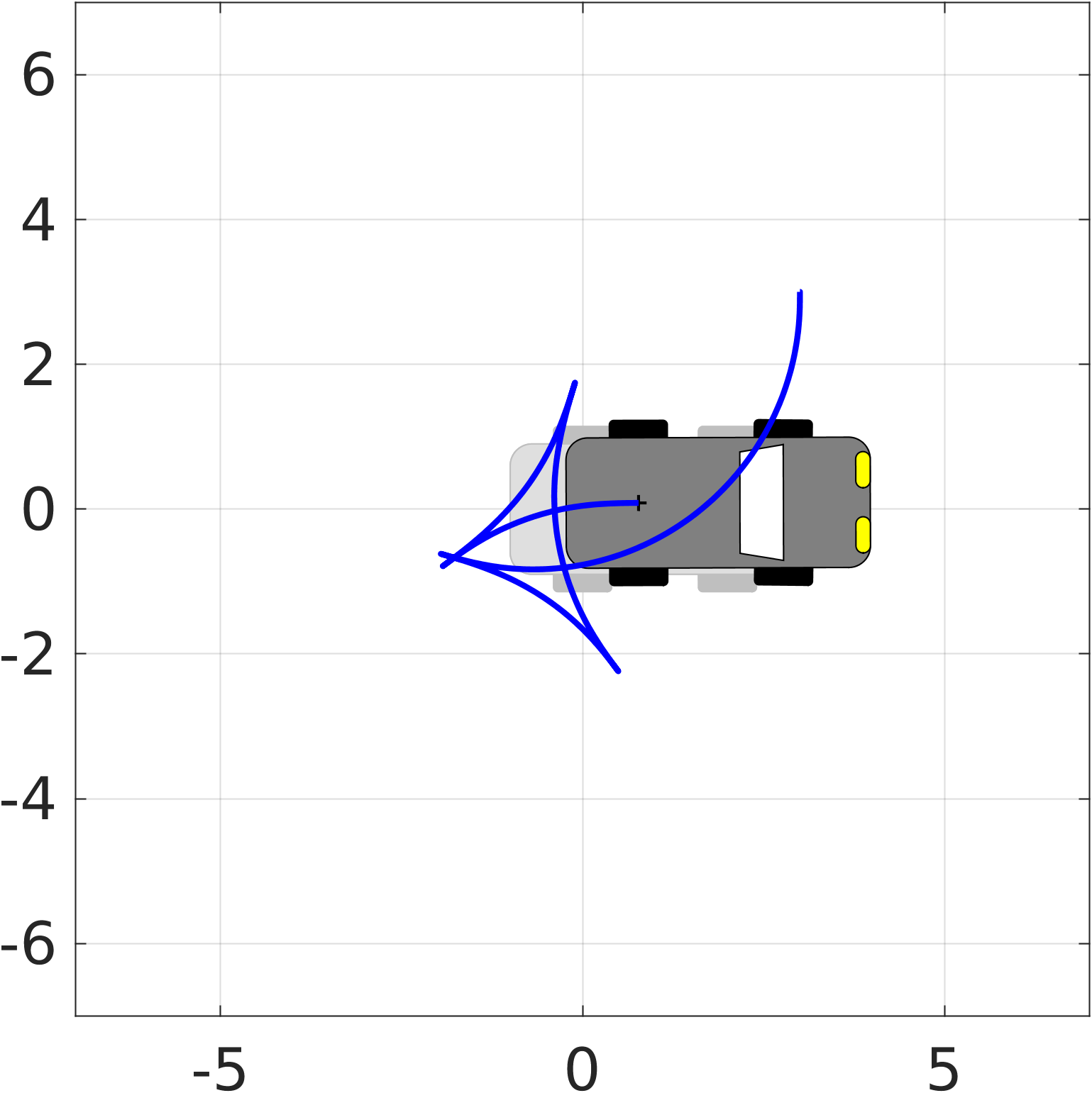}
  }\hfill
  \subfloat[Control input\label{fig:simc}]{
    \includegraphics[width=0.9\linewidth]{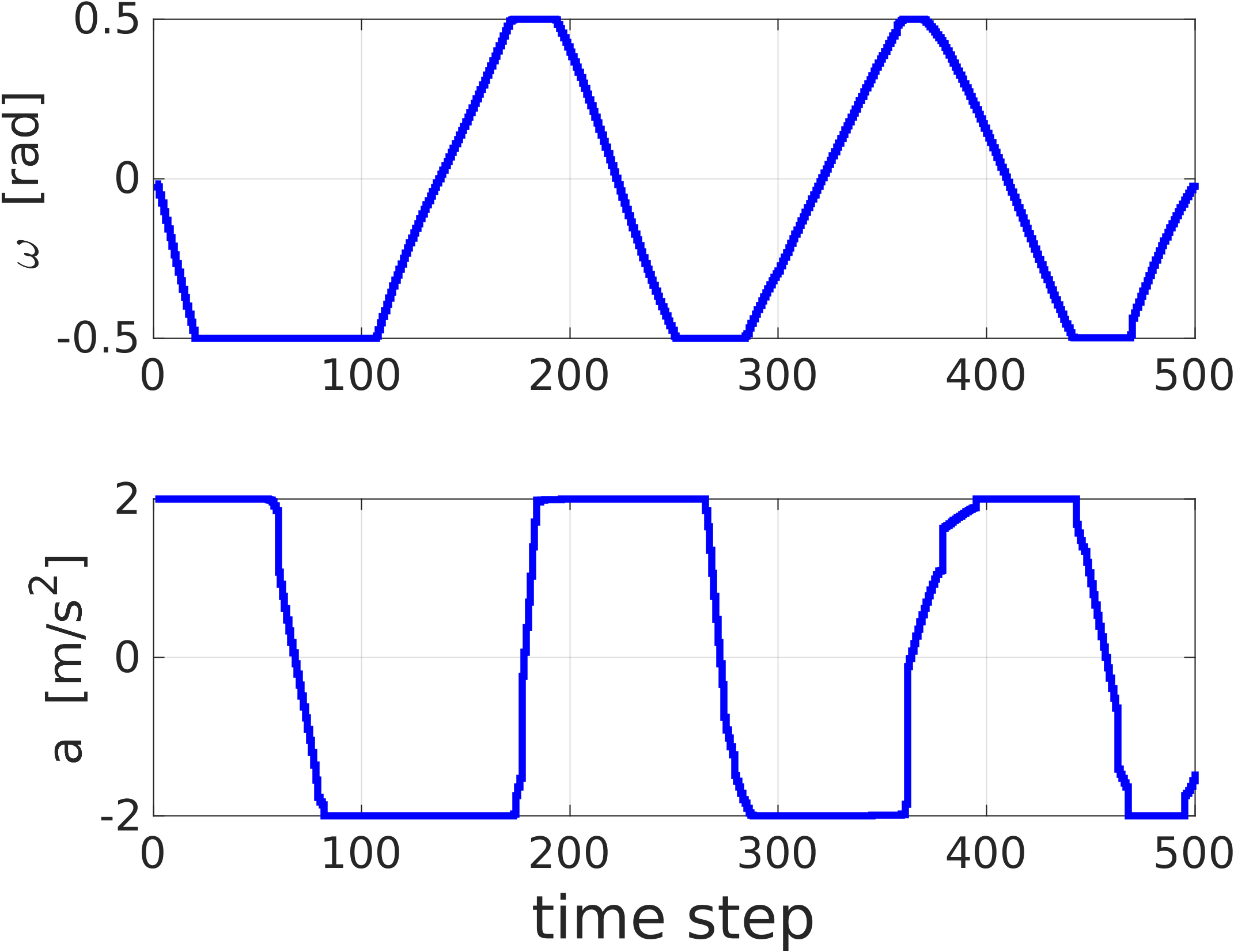}
  }
  \caption{Simulation result for ETS-DDP.}
  \label{fig:simulation}
\end{figure}

\begin{figure}
	\centering
	\includegraphics[width=1\linewidth]{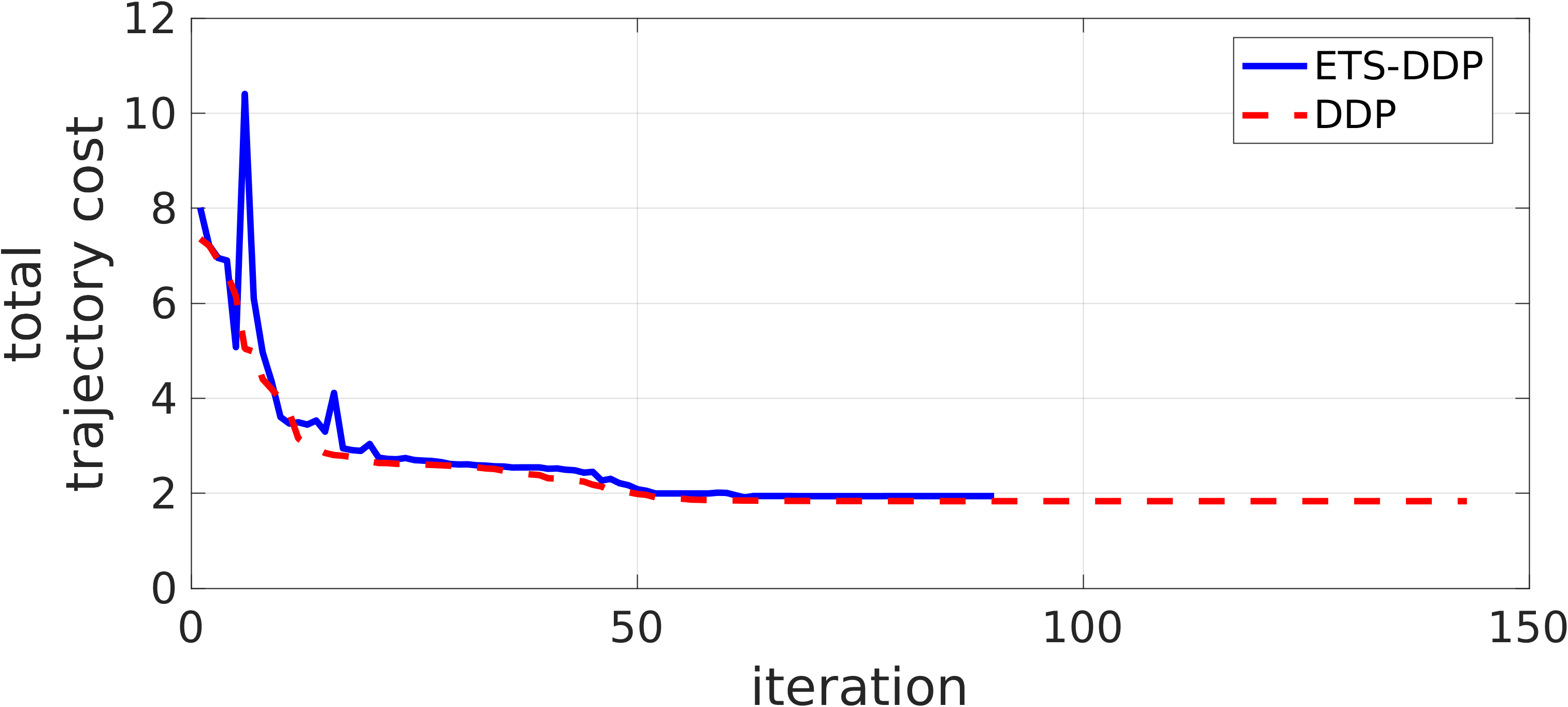}
	\caption{Total cost over iteration associated with conventional DDP (red dashed) and proposed ETS-DDP (blue solid), plotted with the original cost with a target point.}
	\label{fig:costhistories}
\end{figure}

\subsection{Data Collection and Synthesis of the Target Set $\mathcal{C}$}

This subsection explains a detailed procedure for constructing $\mathcal{C}$ in a statistical way as presented in the previous section including data collection.
Before we inferring the parameters required for synthesizing the target set, it is necessary to gather data points $c^{(i)}$, each of which represents a possible (terminal) state of the vehicle with which an expert may conclude that the parking succeeds. We exploit the \textit{rejection sampling} for the data collection conducted through a GUI-type data gathering tool as seen in Fig.~\ref{fig:software}. 
The data gathering tool randomly samples a candidate for the terminal state, say $\widehat{c}$, from the proposal distribution following ${\mathcal{N}}(\widehat{o}, \widehat{\Sigma})$ with mean $\widehat{o}=0$ and variance $\widehat{\Sigma}=\text{diag}(0.1, 0.1, ((15\pi/180))^2, 10^{-5})$. 
A human expert then decides whether to accept the sample or not, considering the available parking area (dashed line in Fig.~\ref{fig:software}) and the corresponding shape of the vehicle (solid line in Fig.~\ref{fig:software}).
In our case, with the help of the GUI tool, we collect $86$ data points $c^{(i)}$, all of which are depicted in Fig.~\ref{fig:software}(c).
With $N=86$ data points, we compute $\bar{c}$ and $S$, and take $r$ as in \eqref{eq:parameter_C_2} with the significance level $\alpha = 0.01$.

\begin{remark}
	 The rejection sampling above is formalized as follows. The proposal distribution follows the normal distribution ${\mathcal{N}}(\widehat{o}, \widehat{\Sigma})$. 
	We also assume that the desired distribution follows a normal distribution for convenience. A typical rejection sampling first draws a sample $\hat{c}$ from the proposal distribution with probability density function (pdf) $f$, which envelops a desired distribution's pdf $g$ when rescaled by $M$. Then, it draws another sample $z$ from uniform distribution $U[0,Mf(\hat{c})]$. If the sample $z$ is in the interval $[0,g(\hat{c})]$ then accept the sample $\hat{c}$. If not, it rejects the sample $\hat{c}$. In our settings, this acceptance and rejection procedure is done by human expert.
	$\hfill\square$
\end{remark}

\subsection{Simulation Result}
Figs.~\ref{fig:simulation} and \ref{fig:costhistories} and Table~\ref{tab:method_comparison} present the simulation results on two DDP algorithms, with $x_0 = (3,3,3\pi/2, 0)$. 
In Fig.~\ref{fig:sima}, the blurry car  depicts the desired target state at $(0,0,0,0)^\top$, while the resulting trajectory generated by the proposed ETS-DDP is drawn in blue solid line. 
In Fig.~\ref{fig:costhistories}, we compare the histories of the total cost from the conventional DDP and the proposed ETS-DDP with the same cost function $\sum_t L(x_t^{(m)}-c,u_t^{(m)})+\phi(x_T^{(m)}-c)$.
It is seen that with a sacrifice of optimality, the proposed ETS-DDP generates an acceptable path much faster than the conventional DDP, whose terminal state enters $\mathcal{C}$ and is close enough to $c=0$ (so that one can conclude that the parking succeeds in an empirical sense).
The overall comparison of iterations, final cost, and computation time are shown in Table~\ref{tab:method_comparison}.

\begin{table}[htbp!]
	\centering
	\caption{Comparison of DDP and ETS-DDP}
	\label{tab:method_comparison}
	\begin{tabular}{l|c|c|c|c}
		\hline
		\textbf{Method} & \textbf{Iter.} & \textbf{Time/Iter.} & \textbf{Total time} & \textbf{Cost} \\
		\hline
		DDP  & 144          & 190 ms        & 27.37 s       &  1.83    \\
		ETS-DDP  & \textbf{91} & 215 ms        & \textbf{19.57 s} & 1.94  \\
		\hline
	\end{tabular}
\end{table}

\section{Conclusion}

In this work, we presented an extended version of the DDP called {\it ETS-DDP}, which solves a class of optimal control problems with the state target chosen as an {\it ellipsoid}. 
The orthogonal projection operator and local smoothing play key roles in reformulating the considered problem with set-valued target into a DDP-friendly form. 
We also provided a statistical way to infer the ellipsoidal target set from pre-collected expert demonstrations. 
It was seen via simulation on autonomous parking that the proposed ETS-DDP computes an acceptable trajectory with which the parking task is achieved, running faster than the conventional DDP while the optimality of the resulting path is lost.
 
It is our belief that the ETS-DDP can be successfully applied to a large class of engineering problems where reducing computational burden is more significant than achieving high-precision control, such as locomotion of legged robots. 
Future works include extension to cases of non-ellipsoidal target sets and/or other statistical distances.


\appendix
We first substitute $y$ with $z:=\Sigma^{-1/2}(y-o)$. Then the optimization problem is represented as
\begin{align*}
    \text{minimize}_z\quad &\|\Sigma^{1/2}z-(x-o)\|^2\\
    \text{subject to}\quad &\|z\|^2 = r^2
\end{align*}
Now, consider this Lagrangian function:
\begin{equation*}
    \mathcal L(z,\lambda) = \|\Sigma^{1/2}z-(x-o)\|^2+\lambda(\|z\|^2 - r^2)
\end{equation*}
The necessary conditions for possible optima are
\begin{equation*}
    \frac{\partial \mathcal L}{\partial z}=0, \quad
    \frac{\partial \mathcal L}{\partial \lambda}=0,
\end{equation*}
leading these normal equations:
\begin{equation*}
    (\Sigma+\lambda I)z=\Sigma^{1/2}(x-o),\quad \|z\|^2=r^2.
\end{equation*}
By solving these normal equations,
the second equation becomes dependent to $\lambda$ and we interpret it as a function of $\lambda$:
\begin{equation*}
    f(\lambda):=\|z\|^2=r^2
\end{equation*}
We call the equation above as the secular equation.
In \cite{gander1978linear}, characterizations of the solution to our problem with regard to the secular equation were given.
\begin{itemize}
    \item Thm 1. the largest $\lambda$ determines solution.
    \item Thm 2. at most one $\lambda > 0$.
\end{itemize}
\cite{gander1978linear} and \cite{reinsch1967smoothing} also covered numerically stable Newton iteration for root finding of secular equations. Instead of using the standard equation for root finding:
\begin{equation*}
    f(\lambda)-r^2=0,
\end{equation*}
the version below showed better global convergence:
\begin{equation*}
    \frac{1}{\sqrt{f(\lambda)}}-\frac{1}{r}=0.
\end{equation*}
Then the update rule of $\lambda_k$ is
\begin{equation*}
    \lambda_{k+1}=\lambda_k - \frac{f(\lambda_k)-r^2}{f'(\lambda_k)}\frac{2\frac{\sqrt{f(\lambda_k)}}{r}}{1+\frac{r}{\sqrt{f(\lambda_k)}}}.
\end{equation*}
The proposed algorithm that numerically solves the orthogonal projection operator is given in Algorithm \ref{alg:sec}.
\begin{algorithm}[htbp!]
\caption{Solving the Secular Equation \cite{gander1978linear}}\label{alg:sec}
\begin{algorithmic}[1]
\Require $\lambda_0 > 0$: Lagrange multiplier, $\Sigma$: dispersion matrix, $o$: the center of ellipsoid, $x$: state obtained from previous iteration
\For{$k=0,1,2,\ldots, K$}
    \State Compute
    \begin{equation*}
        B_k = (\Sigma + \lambda_k I)^{-1}
    \end{equation*}
    \State Compute
    \begin{align*}
        z_k &= B_k \Sigma^{1/2}(x-o)\\
        z_k'&=-B_kz_k\\
        f(\lambda_k)&=z_k^\top z_k\\
        f'(\lambda_k)&=2(z_k')^\top z_k
    \end{align*}
    \State Update $\lambda_{k+1}$
    \begin{equation*}
        \lambda_{k+1}=\lambda_k - \frac{f(\lambda_k)-r^2}{f'(\lambda_k)}\frac{2\frac{\sqrt{f(\lambda_k)}}{r}}{1+\frac{r}{\sqrt{f(\lambda_k)}}}
    \end{equation*}
\EndFor\\
Return $z_K=B_K\Sigma^{1/2}(x-o)$
\end{algorithmic}
\end{algorithm}
After running Algorithm \ref{alg:sec}, the final solution is obtained as
\begin{align*}
    y(x)=\Sigma^{1/2}z_K+o & = \Sigma^{1/2} B_K\Sigma^{1/2}(x-o)+o\\
    & =  \Sigma^{1/2} (\Sigma + \lambda_K I)^{-1}\Sigma^{1/2}(x-o)+o.
\end{align*}
\end{document}